\documentclass[12pt]{article}
\usepackage{amscd, amssymb, latexsym, amsmath}

\newtheorem{pro}{Proposition}

\newenvironment{proof}
{\noindent {\em Proof}} {\hfill $\Box$}

\numberwithin{thm}{section}
\numberwithin{cor}{section}
\numberwithin{pro}{section}
\numberwithin{lem}{section}
\numberwithin{dfn}{section}
\numberwithin{rem}{section}
\numberwithin{equation}{section}

\newcommand{\R}{\mathbb R}
\newcommand{\C}{\mathbb C}

\begin{document}
\title{A Bernstein Type Result for Special \\Lagrangian Submanifolds
}
\author{Mao-Pei Tsui \& Mu-Tao Wang }
%\address{Department of Mathematics\\
%Columbia University\\
%New York , NY 10027}
\date{October 25, 2001, revised Jan 21, 2002}
\maketitle
%\centerline{Preliminary version}

\vskip 10pt \centerline{email:\,tsui@math.columbia.edu,
mtwang@math.columbia.edu}

\begin{abstract}
Let $\Sigma$ be a complete minimal Lagrangian submanifold of
$\C^n$.
 We identify regions in the Grassmannian of Lagrangian
subspaces so that whenever the image of the Gauss map of $\Sigma$
lies in one of these regions, then $\Sigma$ is an affine space.
\end{abstract}

\section{Introduction}

The well-known Bernstein theorem states any complete minimal
surface that can be written as the graph of a function on $\R^2$
must be a plane. This type of result has been generalized in
higher dimension and codimension under various conditions. See
\cite{eh} and the reference therein for the codimension one case
and \cite{b}, \cite{fc}, and \cite{jx} for higher codimension
case. In this note, we prove a Bernstein type result for complete
minimal Lagrangian submanifolds of $\C^n$. We remark that Jost-Xin
\cite{jx2} obtained similar results from a somewhat different
approach.

Recall a submanifold $\Sigma$ of $\C^n$ is called Lagrangian if
the K\"ahler form $\sum_{i=1}^n dx^i\wedge dy^i$ restricts to zero
on $\Sigma$. If $\Sigma$ happens to be the graph of a
vector-valued function from a Lagrangian subspace $L$ to its
complement $L^\perp$ in $\C^n$. Rotating $\C^n$ by a element in
$U(n$), we may assume $L$ is the $x^i$ subspace and $L^\perp$ is
the $y^i$ subspace. In this case, there exists a smooth function $
F : \R^n \to \mathbb R $ such that $\Sigma$ is defined by the
gradient of $F$, $\nabla F$.  The minimal Lagrangian equation can
be written in terms of $F$.

\begin{equation}
\label{eq1.1}
 Im  (\det ( (I + i \mbox{ Hess } (F))) =\textit{ constant}
\end{equation}
where $I$ = identity matrix and Hess $ F = ( \frac{\partial^2 F}{\partial x^i \partial x^j  })$. \\

Such minimal submanifolds were first studied by Harvey and Lawson
\cite{hl} in the context of calibrated geometry. In fact, they are
calibrated by $n$ forms of the type
$Re(e^{i\theta}dz^1\wedge\cdots \wedge dz^n)$ for some constant
$\theta$. They are usually referred as special Lagrangian
submanifold (SLg) in literature in a more general sense.
Recently, Strominger-Yau-Zaslow \cite{syz} established a
conjectural relation of fibrations by special Lagrangian tori with
mirror symmetry.

In terms of (\ref{eq1.1}), a Bernstein type question is to
determine under what conditions an entire solution $F$ becomes a
quadratic polynomial.

The results in this paper imposes conditions on the image of the
Gauss map of $\Sigma$. Recall the set of all Lagrangian subspaces
of $\C^n$ is parametrized by the Lagrangian Grassmannian
$U(n)\slash SO(n)$. The Gauss map of a Lagrangian submanifold
$\gamma:\Sigma\mapsto U(n)\slash SO(n)$ assigns to each $x\in
\Sigma$ the tangent space at $x$, $T_x\Sigma$.

A particular subset of the Lagrangian Grassmannian consists of the
graphs of any symmetric linear transformation from $\R^n$ to
$\R^n$. These can be considered as Lagrangians defined by the
gradient of quadratic polynomials on $\R^n$.

For any $K>0$, let $\mathfrak{B}_K$ to be the subset of the
Lagrangian Grassmannian consisting of graphs of symmetric linear
transformations $L:\R^n\mapsto \R^n$ with eigenvalues
$|\lambda_i|\leq K$ for each $i$. We remark that if the Gauss map
of $\Sigma$ lies in $\mathfrak{B}_K$ then $\Sigma$ is the graph of
$f:\R^n\mapsto \R^n$ with uniformly bounded $|df|$.

\vskip 10pt \noindent {\bf Theorem A} {\it  Denote by $\Xi$ the
subset of the Lagrangian Grassmannian consisting of graphs of
symmetric linear transformations $L:\R^n\mapsto \R^n$ with
eigenvalues $\lambda_i\lambda_j\geq -1$ for any $i, j$. Let
$\Sigma$ be a complete minimal Lagrangian submanifold of $\C^n$.
Suppose there exists an element $g\in U(n)$  such that the image
of the Gauss map of $g(\Sigma)$ lies in $\Xi \cap \mathfrak{B}_K$,
then $\Sigma$ is an affine space. } \vskip 10pt

 We remark that  the
gradient of $g(\Sigma)$ is not necessarily bounded.

Indeed, the most general theorem of this type is the following.

Let $\mathfrak{M}$ be the set of graphs of all symmetric linear
transformation $L:\R^n\mapsto\R^n$ whose eigenvalues $(\lambda_i)$
satisfy the following two conditions:

\begin{enumerate}
\item\[F(h_{ijk})=\sum_{i,j, k} h_{ijk}^2 + \sum_{k, i} \lambda_{
i}^2 h_{iik}^2 +2\sum_{k, i<j}\lambda_{i}\lambda_{j}h_{i j
k}^2\geq 0\] for any trace-free symmetric three tensor  $h_{ijk}$.

\item \[F(h_{ijk})= 0\]if and only if $h_{ijk}=0$ for all $i, j,k$.
\end{enumerate}

 Here $h_{ijk}$ is any element in
$\otimes^3 \R^n$ that is symmetric in $i,j $ and $k$. $h_{ijk}$
being trace-free means $\sum_{i=1}^n h_{iik}=0$ for any $k$. In
fact, $h_{ijk}$ corresponds to the second fundamental form of a
Lagrangian submanifold. The trace-free condition corresponds to
vanishing mean curvature vector. It is clear that $\Xi$ is a
subset of $\mathfrak{M}$.

 \vskip 10pt \noindent {\bf Theorem B} {\it The conclusion for
Theorem A holds for $\mathfrak{M}_K$ , the subset of the
Lagrangian Grassmannian consisting of graphs of symmetric linear
transformations in $\mathfrak{M}\cap \mathfrak{B}_K.$} \vskip 10pt

These theorems are proved by maximum principle. When $\Sigma$ is
the graph over a Lagrangian subspace $L$, we calculate the
Laplacian of  $\ln *\Omega$ where $*\Omega$ is the Jacobian of the
projection from $\Sigma$ to $L$. This is a positive function and
when the Gauss map of $\Sigma$ satisfies the above conditions it
is indeed superharmonic. The parabolic version of this equation
was first derived in \cite{mu3} in the study of higher
co-dimension mean curvature flow.

We wish to thank Professor D. H. Phong and Professor S.-T. Yau for
their  encouragement and  support.

\section{Proof of Theorem}

Let $\Sigma$ be a complete submanifold of $\R^{2n}$.
 Around any point $p\in \Sigma$, we choose orthonormal frames
$\{e_i\}_{i=1\cdots n}$ for $T\Sigma$ and
$\{e_\alpha\}_{\alpha=n+1, \cdots, 2n}$ for $N\Sigma$, the normal
bundle of $\Sigma$ . The convention that $i, j, k, \cdots$ denote
tangent indexes and $\alpha, \beta, \gamma \cdots $ denote normal
indexes is followed.

The second fundamental form of $\Sigma$ is denoted by  $h_{\alpha
ij} =<\nabla_{e_i} e_j, e_\alpha>$.

The following formula was essentially derived in \cite{mu3}. To
apply to the current situation, we note that minimal submanifold
corresponds to stationary phase of mean curvature flow.

\begin{pro}\label{equation}
Let $\Sigma$ be a the graph of $f:\R^n\mapsto \R^n$ and
$(\lambda_i)$ be the eigenvalues of $\sqrt{(df)^T df}$. If
$\Sigma$ is a minimal submanifold,  then $*\Omega=\frac{1}
{\sqrt{\prod_{i=1}^n(1+\lambda_{i}^2)}}$ satisfies the following
equation.
\begin{equation}\label{eq}
\begin{split}
&\Delta *\Omega= -*\Omega\{\sum_{\alpha, i, k} h_{\alpha
ik}^2-2\sum_{k, i<j} \lambda_{ i}\lambda_{ j}h_{n+i,ik} h_{n+j,
jk} +2\sum_{k, i<j}\lambda_{i}\lambda_{j} h_{n+j, ik} h_{n+i,
jk}\}
\end{split}
\end{equation}
where $\Delta $ is the Laplace operator of the induced metric on
$\Sigma$.

\end{pro}

Geometrically, $*\Omega$ is the Jacobian of the projection from
$\Sigma$ to the domain $\R^n$.

\vskip 10pt
\begin{proof}\textit{\,of Theorem A.} First we show if the Gauss map
of $\Sigma$ lies in $\Xi\cap \mathfrak{B}_K$, then $\Sigma$ is an
affine space. The general case follows from the following
observation: if $g\in U(n)$ then $g(\Sigma)$ is again a minimal
Lagrangian submanifold.

 We rewrite equation (\ref{eq}) in the Lagrangian case.  Hence the tangent bundle is
canonically isomorphic to the normal bundle by the complex
structure $J$. We define

 \[h_{ijk}=<\nabla_{e_i} e_j, J(e_k)>\]
then $h_{ijk}$ is symmetric in $i,j$ and $ k$.

The Lagrangian condition also implies $<df(X), J(Y)>$ is symmetric
in $X, Y$. We can find an orthonormal basis $\{e_i\}_{i=1\cdots
n}$ for $T_p\Sigma$ so that $df(e_i)=\lambda_i J(e_i )$ and
$\{J(e_i)\}_{i=1\cdots n}$ becomes an orthonormal basis for the
normal bundle. Equation (\ref{eq}) becomes

\begin{equation}\label{eq1}
\begin{split}
&\Delta *\Omega= -*\Omega\{\sum_{i,j, k} h_{ijk}^2-2\sum_{k, i<j}
\lambda_{ i}\lambda_{ j}h_{iik} h_{j jk} +2\sum_{k,
i<j}\lambda_{i}\lambda_{j} h_{j ik} h_{i jk}\}
\end{split}
\end{equation}

We shall calculate
\begin{equation}\label{eq2}
\Delta ( \ln *\Omega ) = \frac{ *\Omega \Delta ( *\Omega ) -
|\nabla *\Omega|^2}{ |*\Omega|^2 }
\end{equation}

The covariant derivative of $*\Omega$ can be calculated as in
equation (3.1) in \cite{mu3}.

 \[(*\Omega)_{k} = -*\Omega
(\sum_{i} \lambda_{i}h_{iik} ) \]

Plug this and equation (\ref{eq1}) into equation (\ref{eq2}) and
we obtain

\begin{equation}\label{same}
\Delta ( \ln *\Omega )= - \{\sum_{ i,j, k} h_{ijk}^2 + \sum_{k, i}
\lambda_{ i}^2 h_{iik}^2 +2\sum_{k, i<j}\lambda_{i}\lambda_{j} h_{
ijk}^2\ \}
\end{equation}

If the Gauss map of $\Sigma$ lies in $\Xi$, then it is obvious
that $\Delta ( \ln *\Omega )\leq 0$. The condition
$|\lambda_i|\leq K$ means $\Sigma$ is the graph of a vector-valued
function with bounded gradient. $\sum_{ i,j, k} h_{ijk}^2 +
\sum_{k, i} \lambda_{ i}^2 h_{iik}^2 +2\sum_{k,
i<j}\lambda_{i}\lambda_{j} h_{ ijk}^2=0$ forces $h_{ijk}=0$ for
any $i, j,k$ by symmetry consideration. This immediately implies
any minimal Lagrangian cone satisfies the assumption of the
theorem is flat by maximum principle. For the general case, we can
apply the standard blow-down and Allard regularity theorem to
conclude $\Sigma$ is totally geodesic and thus an affine space.

If $\Sigma$ is minimal Lagrangian, so is any $g(\Sigma)$ for $g\in
U(n)$. This is because $U(n)$ is contained in the isometry group
of $\C^n$ and it preserves the standard K\"ahler form. This
completes the proof of Theorem A.

\end{proof}

\vskip 10pt
\begin{proof}\textit{\,of Theorem B.}

This follows immediately from the definition of the set $\Xi'_K$
and equation (\ref{same}). Because $\Sigma$ is minimal, we only
need to consider trace-free $h_{ijk}$.

\end{proof}

In the rest of the paper, we identify another region of the
Lagrangian Grassmannian where the Bernstein-type theorem also
applies. This is not as general as the region in Theorem B.
However, we expect it will provide a better estimate in future
application.

 \vskip 10pt \noindent {\bf Theorem
C} {\it The conclusion for Theorem A holds for $\Xi'\cap
\mathfrak{B}_K$ where $\Xi'$ is the subset of Lagrangian
Grassmannian consisting of graphs of symmetric linear
transformations $L:\R^n\mapsto \R^n$ with eigenvalues
$\lambda_i\lambda_j+\lambda_i\lambda_k+\lambda_j\lambda_k \geq 0$
for any pairwise distinct $i, j, k$.} \vskip 10pt

\vskip 10pt
\begin{proof}\textit{\,of Theorem C.}

We rewrite the right hand side of equation (\ref{same}).
\begin{equation}\label{laplno}
\begin{split}
& - \{\sum_{ i,j, k} h_{ijk}^2 + \sum_{k, i} \lambda_{ i}^2
h_{iik}^2 +2\sum_{k, i<j}\lambda_{i}\lambda_{j}h_{i j k}^2\ \}
\\= & - \{\sum_{ i,j, k}
h_{ ijk}^2  + \sum_{ i}\lambda_{ i}^2 h_{iii}^2 + \sum_{ i < k}(
\lambda_{ i}^2 + 2\lambda_{ i}\lambda_{ k} )h_{iik}^2 + \sum_{ i
> k} ( \lambda_{ i}^2 + 2\lambda_{ i}\lambda_{ k} ) h_{iik}^2
\\& +2\sum_{ i<j < k }( \lambda_{i}\lambda_{j} + \lambda_{ j}\lambda_{ k} + \lambda_{ k}\lambda_{ i} ))h_{i j k}^2 \}
\end{split}
\end{equation}

Since $\Sigma$ is minimal, the mean curvature vector $\sum_{i=1}^n
h_{iik}=0$ for each $k$ , we have
\begin{equation}\label{H}
\begin{split}
&\sum_{ i}\lambda_{ i}^2 h_{iii}^2
\\ = &\sum_{ i < j }\lambda_{ i}^2 h_{ijj}^2 +\sum_{
i > j }\lambda_{ i}^2 h_{ijj}^2 + 2 \sum_{ i \neq j, i \neq l, j <
l }\lambda_{ i}^2 h_{ijj} h_{ill}
\\= &\sum_{ k < i}\lambda_{ k}^2 h_{kii}^2 +\sum_{
k > i }\lambda_{k}^2 h_{kii}^2 + 2 \sum_{ i \neq j, i \neq l, j <
l }\lambda_{ i}^2 h_{ijj} h_{ill}
\end{split}
\end{equation}

Plug equations (\ref{H}) into (\ref{laplno}),

\begin{equation}
\begin{split}
&\Delta ( \ln *\Omega )
\\= & - \{\sum_{i, j, k}
h_{ ijk}^2 + 2 \sum_{  i \neq j ,i \neq l, j < l }\lambda_{ i}^2
h_{ijj} h_{ill} + \sum_{ i < k}( \lambda_{ i}^2 + 2\lambda_{
i}\lambda_{ k} + \lambda_{ k}^2)h_{iik}^2 +
\\ &\sum_{ i > k} ( \lambda_{ i}^2 + 2\lambda_{ i}\lambda_{ k} + \lambda_{ k}^2) h_{iik}^2
 +2\sum_{ i<j < k }( \lambda_{i}\lambda_{j} + \lambda_{ j} \lambda_{ k} + \lambda_{ k}\lambda_{ i} ))h_{i j k}^2 \}
\\ =  & - \{\sum_{ i, j,k}
h_{ ijk}^2 + 2 \sum_{ i \neq j ,i \neq l, j < l }\lambda_{ i}^2
h_{ijj} h_{ill} + \sum_{   p \neq q}( \lambda_{ p} +
\lambda_{q})^2 h_{pqq}^2 +
\\ &
 +2\sum_{ i<j < k }( \lambda_{i}\lambda_{j} + \lambda_{ j} \lambda_{ k} + \lambda_{ k}\lambda_{ i} )h_{i j k}^2\
\}
\end{split}
\end{equation}

\begin{equation}
\begin{split}
& 2 \sum_{ i \neq j ,i \neq l, j < l }\lambda_{ i}^2 h_{ijj}
h_{ill} + \sum_{   p \neq q}( \lambda_{ p} + \lambda_{q})^2
h_{pqq}^2
\\ & \ge  \sum_{ i \neq j ,i \neq l, j < l }2\lambda_{
i}^2 h_{ijj} h_{ill} + ( \lambda_{ i} + \lambda_{ l})^2 h_{i l
l}^2 + ( \lambda_{ i} + \lambda_{j})^2 h_{ijj}
\end{split}
\end{equation}

If  $ j \neq i  \neq l$ and $ \lambda_{i}\lambda_{j} + \lambda_{
j} \lambda_{ l} + \lambda_{ l}\lambda_{ i} \ge 0 $ then

\[ 2\lambda_{ i}^2 h_{ijj} h_{ill} + ( \lambda_{ i} + \lambda_{
l})^2 h_{i l l}^2 + ( \lambda_{ i} + \lambda_{j})^2 h_{ijj} \ge
0.\]

Thus $\Delta ( \ln *\Omega ) \le -|A|^2 $. The rest is identical
to that of Theorem A and B.

\end{proof}

We remark the condition
$\lambda_i\lambda_j+\lambda_i\lambda_k+\lambda_j\lambda_k \geq 0$
for any pairwise distinct $i, j, k$ is void in two-dimension.
Indeed, this is true even without the Lagrangian assumption and
hence rediscover the results of \cite{b} (see also \cite{fc}) that
an non-parametric minimal cone of dimension three must be flat.

\footnote{After this paper was finished, we were informed by Yu
Yuan that he has also derived formula (2.4) from a different point
of view. In his paper "A Bernstein problem for special Lagrangian
equations", Yu Yuan had the following interesting observation: the
linear transformation $(x^i, y^i)\mapsto
(\frac{x^i+y^i}{\sqrt{2}}, \frac{-x^i+y^i}{\sqrt{2}})$ takes a
convex function $F$ to a function $\overline{F}$  with $-I\leq D^2
\overline{F} \leq I$. Since this transformation (so called Lewy
transformation) is an element of $U(n)$, our Theorem A implies a
convex entire solution to equation (1.1) is a quadratic
polynomial. This result was also proved in Yu Yuan's paper. We
would like to thank Yu Yuan for sending us his preprint before
publication .
 }

\end{document}